\documentclass{amsart}
\usepackage{amscd}
\newtheorem{theorem}{Theorem}[section]

\newtheorem{corollary}[theorem]{Corollary}
\newtheorem{proposition}{Proposition}
\newtheorem{definition}{Definition}
\newtheorem{remark}{Remark}

\begin{document}

\begin{center}
{\huge Degeneration, rigidity and irreducible components of Hopf algebras } \\

\bigskip

\bigskip

Abdenacer MAKHLOUF \\

\medskip

Laboratoire de Math\'{e}matiques et Applications\\
 4, rue des
fr\`{e}res Lumi\`{e}re 68093 Mulhouse, France

Email : N.Makhlouf@univ-mulhouse.fr
\end{center}
\bigskip

\bigskip

\textit{The aim of this work is to discuss the concepts of degeneration, deformation and
rigidity of Hopf algebras and to apply them to the geometric study of
the varieties of Hopf algebras. The main result is the description of the }$n$-
\textit{dimensional
rigid Hopf 
algebras and the irreducible components for }$n=p^{2},$ $p$ \textit{is prime number, and
for } 
$n<14.$\textit{ This paper is organized as follows : first
we present the basic concepts of degeneration, deformation and rigidity with examples and some useful
properties.  The second part of this work
is devoted to the geometric description of the algebraic varieties of
n-dimensional Hopf algebras . In the last section, we give a necessary and sufficient 
condition for the existence of a degeneration of a given Hopf algebra. }\newline

2000 AMS Classification : 16W30, 14D15, 13D10.


\section{Introduction}

The discovery of quantum groups gives impulse to a strong development in the
theory of Hopf algebras. Many articles were devoted to the class of
semisimple Hopf algebras; several characterizations were given and some
fundamental results were established. See for example the surveys \cite{Andrus1} \cite
{Montgomery}. The complete classification of all Hopf algebras
of a fixed dimension $n$ is known only for $n=p$ ($p$ prime) \cite{Zhu1}, $%
n=p^{2}$ ($p$ prime)\cite{Ng} and for certain small dimension $n$, $n<14$, $%
n=15,21,35$. Besides this, substantial results are known for some classes,
like pointed Hopf algebras (see \cite{Andrus-Schneider2}) , and triangular
Hopf algebras (see \cite{Etingof-Gelaki})

This paper is focused on the geometric and local study of the variety of
Hopf algebras of fixed dimension $n$. The only known antecedent is the 
paper \cite{Stefan2}, where D. Stefan considers, in a
cohomological way, the invariant of the irreducible components of semisimple
and cosemisimple Hopf algebras. Here, we approach the varieties of Hopf
algebras, in a more geometric way, via the algebraic deformation theory introduced by
M. Gerstenhaber \cite{Gerstehaber92} and degeneration theory.

The concept of degeneration appeared first in the physics literature. The
question was to show in which sense a group can be a limiting case of other
groups. Degenerations, called also contractions or specialisations, were
introduced for Lie groups by Segal, Inonu and Wigner (1953) \cite
{Inonu-Wigner}. They showed that the Galilei group of classical mechanic is
a limiting case of the Lorentz group corresponding to relativistic
mechanic. Later, Saletan (1960) \cite{Saletan} generalized the notion and
stated a general condition for the existence of degenerations (or
contractions) of Lie algebras.

The notions of deformation and degeneration were used by several autors in
the studies of associative algebras varieties or Lie algebras varieties
(Ancochea-Bermudez, Carles, Gabriel, Gerstenhaber, Goze, Happel, Mazzola, 
Makhlouf, Schaps, Schack...). It was
also used, in the theory of quantum groups, by Celegheni, Giachetti, Sorace
and Tarlini \cite{Celeghini} to define Heisenberg and Euclidean quantum
groups.

We now discuss the contents of this paper. We first recall 
the algebraic variety structure of the set of bialgebras $Bialg_{n}$ and
Hopf algebras $Hopf_{n}$. Then, we describe the linear group action on these
varieties. The section 3 is devoted to the degeneration of Hopf algebras, 
we give some examples and show that the graded Hopf algebra associated to a filtred Hopf algebra
 is always a degeneration of this Hopf algebra. 
  The section 4 is devoted to deformation and rigidity of 
Hopf algebras, we give some useful properties and prove that if a finite dimensional 
Hopf algebra is rigid then
its dual is rigid. In the section 5, we study the deformations and
degenerations of Hopf algebras appearing in known classification. We consider first the
algebraic varieties $
Hopf_{p^{2}}$ where $p$ is any prime number. We prove that every  Hopf algebra, in these
algebraic varieties, is rigid.
Therefore these algebraic varieties are union of Zariski open orbits. And next 
we  show that, for $%
\mathit{n}<14,$ every n-dimensional Hopf algebra is rigid except  $A_{c_{4}}^{^{ \prime }}$ in $Hopf_{8}$
and $A_{0}$ , $B_{1}$ in $Hopf_{12}$. We prove that $A_{c_{4}}^{^{\prime }}
$ is a degeneration of $A_{c_{4}}^{^{\prime \prime }}$ and $A_{0}$ is a
degeneration of $A_{1}$ and by duality $A^{\star}_{0}\simeq B_{1}$ is a degeneration of
$A^{\star}_{1}$. Then, we describe the irreducible components of $
Hopf_{n}$ for $\mathit{n}<14$.
In the last section, we give a necessary and sufficient condition for the existence of a
degeneration of a given Hopf algebra.

\textbf{Acknowledgement. }

I wish to express my thanks to N. Andruskiewitsch for the helpful discussions
and his careful reading the original manuscript.

\section{Structure of algebraic varieties}

Throughout this paper $K$ will be an algebraically closed field of
characteristic $0$ and $V$ be an $n$-dimensional vector space over the field 
$K$. Let $H=(V,\mu ,\eta ,\Delta ,\varepsilon ,S)$ be a finite dimensional
Hopf algebra, where $\mu :V\otimes V\rightarrow V$ is the multiplication, $\eta
:K\rightarrow V$ is the unity, $\Delta :V\rightarrow V\otimes V$ is the comultiplication,
$\varepsilon :V\rightarrow k$ is the counity and the endomorphism $S$ is the
antipode. We refer to \cite{MontgomeryL} for the definitions.

Setting a basis \{$e_1,...,e_n\}$ of $V$ where $e_1=\eta \left( 1\right) ,$
we identify the multiplication $\mu $ and the comultiplication $\Delta $
with their $n^3$ \textit{structure constants }$C_{ij}^k$ and $D_i^{jk}\in K,$
where $\mu \left( e_i\otimes e_j\right) =\sum_{k=1}^nC_{ij}^ke_k$ and $\Delta
\left( e_i\right) =\sum_{j=1}^n\sum_{k=1}^nD_i^{jk}e_j\otimes e_k.$ The
counity $\varepsilon $ is identified with its $n$ structure constants $\xi _i,$
where $\varepsilon \left( e_i\right) =\xi _i$. The collection $\left(
C_{ij}^k,D_i^{jk},\xi _i:i,j,k=1,\cdots ,n\right) $ represents a bialgebra
if it satisfies the following polynomial equations : 
$$
(1)\quad \left\{ 
\begin{array}{l}
\sum_{l=1}^nC_{ij}^lC_{lk}^s-C_{il}^sC_{jk}^l=0 \\ 
C_{1i}^j=C_{i1}^j=\delta _{ij}\text{ the Kronecker symbol}
\end{array}
\right. \quad i,j,k,s\in \left\{ 1,..,n\right\}
$$
$$
(2)\quad \left\{ 
\begin{array}{l}
\sum_{l=1}^nD_s^{lk}D_l^{ij}-D_s^{il}D_l^{jk}=0 \\ 
\sum_{l=1}^nD_i^{il}\xi _l=\sum_{l=1}^nD_i^{li}\xi _l=1 \\ 
\sum_{l=1}^nD_i^{jl}\xi _l=\sum_{l=1}^nD_i^{lj}\xi _l=0\ i\neq j
\end{array}
\right. \quad i,j,k,s\in \left\{ 1,..,n\right\}
$$
$$
(3)\quad \left\{ 
\begin{array}{l}
\sum_{l=1}^nC_{ij}^l-\sum_{r,t,p,q=1}^nD_i^{rt}D_j^{pq}C_{rp}^kC_{tq}^s=0 \\ 
D_1^{11}=1,\quad D_1^{ij}=0\ \left( i,j\right) \neq \left( 1,1\right) \\ 
\xi _1=1,\quad \sum_{l=1}^nC_{ij}^l\xi _l=\xi _i\xi _j
\end{array}
\right. \quad i,j,k,s\in \left\{ 1,..,n\right\}
$$

The polynomial relations $\left( 1\right) \left( 2\right) \left( 3\right) $
endow the set of the $n$-dimensional bialgebras, denoted by $Bialg _{n}$,
with a structure of an algebraic variety embedded in $K^{2n^{3}+n}$ which we
do consider here together with its natural structure of an algebraic variety
over $K$.

As a subset of $K^{2n^{3}+n}$, $Bialg_{n}$ may be provided with the Zariski
topology. If $K$ is the complex field the set $Bialg_{n}$ may be
provided with the metric topology of $K^{2n^{3}+n}$ . Recall that the metric
topology is finer than the Zariski topology.

\begin{remark}
We can also fix the counity as in Stephan's work but it seems more general
keeping it free. The set where the counity is fixed is an open subset in $%
Bialg_{n}$ (for metric or Zariski topologies).
\end{remark}

The elements of Bialg$_{n}$ with antipodes define the set of Hopf algebras,
denoted by $Hopf_{n}.$ If $S=\left( S_{ij}\right) $ defines the
antipode,with respect to the basis \{$e_{1},...,e_{n}\}$ of $V,$ then, in addition to
relations $\left( 1\right) \left( 2\right) \left( 3\right) $, we have
the following equations : 
$$
(4)\left\{ 
\begin{array}{c}
\sum_{j,k,r=1}^{n}D_{i}^{jk}S_{rj}C_{rk}^{1}=%
\sum_{j,k,r=1}^{n}D_{i}^{jk}S_{rk}C_{jr}^{1}=\xi _{i} \\ 
\sum_{j,k,r=1}^{n}D_{i}^{jk}S_{rj}C_{rk}^{t}=%
\sum_{j,k,r=1}^{n}D_{i}^{jk}S_{rk}C_{jr}^{t}=0
\end{array}
\right. \quad i\in \left\{ 1,..,n\right\} ,t\in \left\{ 2,..,n\right\}
$$

The set $Hopf_{n}$ is a Zariski open subset of $Bialg_{n}.$

\begin{remark}
The sets $Bialg_{n}$ and $Hopf_{n}$ carry also a structure of scheme.
\end{remark}

\subsection{The GL$_{\text{n}}\left( \text{K}\right) $ action,
orbits}

In the following we will define $GL_{n}\left( K\right) $ action
for Hopf algebras. Naturally, we have a similar definition for bialgebras.

Geometrically, a point $\left( C_{ij}^{k},D_{i}^{jk},\xi _{i}:i,j,k=1,\cdots
,n\right) $ of $K^{2n^{3}+n}$ satisfying $\left( 1\right) $ $\left( 2\right) 
$ $\left( 3\right) $ $\left( 4\right) $ where the matrix $\left(
S_{ij}\right) $ defines the antipode, represents an $n$-dimensional Hopf
algebra $H$, along with a particular choice of basis. A change of basis in $%
H $ may give rise to a different point of $Hopf_{n}$. Let $H=(V,\mu ,\eta
,\Delta ,\varepsilon ,S)$ be a Hopf algebra, The ''structure transport''
action is defined on $H$ by the following action of $GL_{n}\left( K\right) $ 
$$
\begin{array}{c}
GL_{n}\left( K\right) \times Hopf_{n}\longrightarrow Hopf_{n} \\ 
\left( f,H\right) \longrightarrow f\cdot H
\end{array}
$$
$\forall X,Y\in V$

$$
\begin{array}{c}
\left( f\cdot \mu \right) \left( X\otimes Y\right) =f^{-1}\left( \mu \left(
f\left( X\right) \otimes f\left( Y\right) \right) \right) \\ 
(f\cdot \Delta )\left( X\right) =f^{-1}\otimes f^{-1}\left( \Delta \left(
f\left( X\right) \right) \right) \\ 
\left( f\cdot \varepsilon \right) \left( X\right) =\varepsilon \left(
f\left( X\right) \right)
\end{array}
$$

The \textit{orbit} of a Hopf algebra $H$ is given by $\vartheta \left(
H\right) =\left\{ f\cdot H,\quad f\in GL_n\left( K\right)
\right\} .$ The orbits are in 1-1-correspondence with the isomorphism
classes of $n$-dimensional Hopf algebras. 
The stabilizer subgroup of $H$ $\left( stab\left( H\right) =\left\{ f\in
GL_n\left( K\right) :H=f\cdot H\right\} \right) $ is exactly $Aut\left(
H\right) ,$ the automorphism group of $H$.
The orbit $\vartheta \left( H\right) $ is identified with the homogeneous space $GL_n\left(
K\right) /Aut\left( H\right) $. Then 
$$
\dim \vartheta \left( H\right) =n^2-\dim Aut\left( H\right)
$$

The orbit $\vartheta \left( H\right) $ is provided, when $K=\mathbb{C}$ (a complex
field), with the structure of a differentiable manifold. In fact, $\vartheta
\left( H\right) $ is the image through the action of the Lie group $%
GL_{n}\left( K\right) $ of the point $H$, considered as a point of $%
Hom\left( V\otimes V,V\right) \times Hom\left( V,V\otimes V\right) $.

The Zariski open orbits have a special interest in the geometric study of $%
Hopf_{n}$. It corresponds to a so called rigid Hopf algebras. The orbit's
closure of a rigid Hopf algebra determines an irreducible component of $%
Hopf_{n}$. Many but not all components of $Hopf_{n}$ are orbits closures of
rigid Hopf algebras, for example there are infinitely many isomorphism
classes for $\dim H=p^{4}\ \left( p\text{ odd and prime}\right) $\cite
{Beattie}\cite{Andrus-Schneider}\cite{Gelaki}. This parameter family cannot be in the
closure of a rigid Hopf algebra. 
\section{Degenerations}

\subsection{Definition}

\textit{Let }$H_{0}\ $\textit{and }$H$\textit{\ be two n-dimensional Hopf
algebra. We say \ that }$H_{0}$\textit{\ is a\ degeneration\ \ of\ }$H$%
\textit{\ if }$H_{0}$\textit{\ is in }$\overline{\vartheta \left( H\right) }$%
\textit{\ , the Zariski closure of the orbit of }$H.$
\begin{remark}
\begin{itemize}
\item  Let $t$\ be a parameter in $K$\ and $\{f_{t}\}_{t\neq 0}$\ be a
family of continuous invertible linear maps on $V$\ over $K$\ and $H=\left(
V,\mu ,\Delta ,\eta ,\varepsilon ,S\right) $\ be a Hopf algebra over $K.$\
The limit (when it exists) of a sequence $f_{t}\cdot H$\ , $%
H_{0}=\lim_{t\rightarrow 0}f_{t}\cdot H$ , is a degeneration of $H$ in the
sense that $H_{0}$\ is in the Zariski closure of the set $\left\{ f_{t}\cdot
H\right\} _{t\neq 0}.$
\\  The multiplication and the comultiplication $\mu _{0}$ and $\Delta
_{0}$ of $H_{0}$ are given by 
$$
\begin{array}{c}
\mu _{0}=\lim_{t\rightarrow 0}f_{t}\cdot \mu =\lim_{t\rightarrow
0}f_{t}^{-1}\circ \mu \circ f_{t}\otimes f_{t} \\ 
\Delta _{0}=\lim_{t\rightarrow 0}f_{t}\cdot \Delta =\lim_{t\rightarrow
0}f_{t}^{-1}\otimes f_{t}^{-1}\circ \Delta \circ f_{t}
\end{array}
$$
\item The multiplication $\mu _{t}=f_{t}^{-1}\circ \mu \circ f_{t}\otimes f_{t}$
and the comultiplication $\Delta _{t}=f_{t}^{-1}\otimes f_{t}^{-1}\circ
\Delta \circ f_{t}$ satisfy the conditions of bialgebra, then when $t$
tends to $0$ the conditions remain satisfied.

\item  In the point 1, the linear map $f_{t}$\ is invertible when $t\neq 0$ \ and may be
singular when $t=0.$ Then, one may obtain by degeneration a new Hopf algebra.

\item  Geometrically, $H_{0}$ is a degeneration of $H$ means that $H_{0}$ and $H$ belong to the same
irreducible component.

\item  When $K$ is the complex field, the multiplication and
the comultiplication given by the limit, follows from a limit of the structure constants,
using the metric topology. In fact, $f_{t}\cdot \mu $ and $f_{t}\cdot \Delta $
correspond to a change of basis when $t\neq 0$. when $t=0$, they give eventually a new
point in $Hopf_{n}\subset K^{2n^{3}+n}$.

\item  The same definitions and remarks hold for bialgebras.
\end{itemize}
\end{remark}
\subsection{Example}

In $Hopf_{8}$, the Hopf algebra $H_{0}$ defined by $\frac{K\left\langle
x,g\right\rangle }{\left\langle g^{4}-1,x^{2},gx+xg\right\rangle },\quad
\Delta (g)=g\otimes g,\quad \Delta (x)=x\otimes g+1\otimes x,\quad \varepsilon
(x)=0,\quad \varepsilon (g)=1,\quad S(x)=-xg^{3},S(g)=g^{3}$, is a degeneration of
the Hopf algebra $H_{1}$ defined by $\frac{K\left\langle x,g\right\rangle }{%
\left\langle g^{4}-1,x^{2}-g^{2}+1,gx+xg\right\rangle },\quad $ with the same
counity, antipode and coalgebra structure. In fact, the family $H_{t}$
defined by $\frac{K\left\langle x,g\right\rangle }{\left\langle
g^{4}-1,x^{2}+t(1-g^{2}),gx+xg\right\rangle }$ where $t$ is a parameter, is
isomorphic to $H_{1}$ when $t\neq 0$ and tends to $H_{0}$ when $t$ tends to $%
0$.
\subsection{The graded Hopf algebra}
\begin{theorem}

Let \( H \) be a Hopf K-algebra and \( H_{0} \)\( \subseteq H_{1}\subseteq \cdots \subseteq H_{n}\subseteq \cdots  \)
a Hopf algebra filtration of \( H \), \( H=\cup _{n\geq 0}H_{n} \)
and \( \Delta H_{n}=\sum _{m\geq 0}H_{m}\otimes H_{n-m} \). Then the
graded Hopf algebra \( gr(H)=\oplus _{n\geq 1}H_{n}/H_{n-1} \) is
a degenration of \( H \).
\end{theorem}
\begin{proof}

Let \( H[[t]] \) be a power serie ring in one variable $t\in K$ over \( H \), 
\( H[[t]]=H\otimes K[[t]]=\oplus _{n\geq 0}H\otimes t^{n} \).

We denote by \( H_{t} \) the Rees algebra associated to the filtred
Hopf algebra \( H \), \( H_{t}=\sum _{n\geq 0}H_{n}\otimes t^{n} \).
The Rees algebra \( H_{t} \) is contained in the algebra \( H[[t]] \).

For every \( \lambda \in K \), we set \( H_{(\lambda )}=H_{t}/((t-\lambda )\cdot H_{t}) \).
For \( \lambda =0 \), \( H_{(0)}=H_{t}/(t\cdot H_{t}) \). The Hopf algebra 
\( H_{(0)} \) corresponds to the graded algebra $gr(H)$ and $H_{(1)}$ is
isomorphic to $H$. In fact, we suppose that the parameter
\( t \) commutes with the elements of \( H \) then
 \( t\cdot H_{t}= \)\( \oplus _{n\geq 0}H_{n}\otimes t^{n+1} \). 
It follows that $
H_{(0)}=H_{t}/(t\cdot H_{t})=
(\sum _{n}H_{n}\otimes t^{n})/(\oplus _{n}H_{n}\otimes t^{n+1})=
\oplus _{n\geq 0}H_{n}/H_{n-1}=gr(H)$.
By using the linear map from \( H_{t} \) to \( H \) where the image
of \( a_{n}\otimes t^{n} \) is \( a_{n} \), we have also 

\( H_{(1)}=H_{t}/((t-1)\cdot H_{t})\cong H \).

If \( \lambda \neq 0 \),
the change of parameter \( t=\lambda T \) shows that \( H_{(\lambda )} \)
is isomorphic to \( H_{(1)} \). This ends the proof that \( H_{(0)}=gr(H) \)
is a degeneration of \( H_{(1)}\cong H \). 
\end{proof}

\section{Deformation and rigidity of Hopf algebras}
In this section, we recall the algebraic deformation and the rigidity notions 
introduced by Gerstenhaber. We show the connection between degeneration and deformation
and give some remarks and properties useful in the geometric
study of $Hopf_{n}$. We prove in particular that the rigidity of a finite dimensional Hopf algebra is
equivalent to the rigidity of its dual.
\subsection{Algebraic deformation} 
The notion of deformation is in some sense the dual notion of the
degeneration. Let $H=\left( V,\mu _0,\Delta _0,\eta _0,\varepsilon _0,S_0\right) $ be a
Hopf algebra over a field $K$. Let $K\left[ \left[ t\right] \right] $ be the
power series ring in one variable $t$. Let $V\left[ \left[ t\right] \right] $
be the extension of $V$ by extending the coefficient domain from $K$ to $%
K\left[ \left[ t\right] \right] $. Then $V\left[ \left[ t\right] \right] $
is a $K\left[ \left[ t\right] \right] $-module and $V\left[ \left[ t\right]
\right] =V\otimes _KK\left[ \left[ t\right] \right] $.  A deformation of $H$ is a one parameter family $H_{t}=\left( V\left[ \left[
t\right] \right] ,\mu _{t},\Delta _{t},\eta _{t},\varepsilon
_{t},S_{t}\right) $. Since the unit, counit and the antipode are preserved
by deformation \cite{Gerstehaber92}, it follows that a deformation of $%
H=\left( V,\mu _{0},\Delta _{0},\eta _{0},\varepsilon _{0},S_{0}\right) $
can be considered as a pair of deformations $\left( \mu _{t},\Delta
_{t}\right) $ which together give on $V\left[ \left[ t\right] \right] $ the
structure of bialgebra over $K\left[ \left[ t\right] \right] $. By $K\left[ \left[ t\right] \right] $-linearity the morphisms $\mu _t,\Delta
_t$ are determined by their restrictions to $V\otimes V$ :

$$
\mu _t: 
\begin{array}{l}
V\otimes V\rightarrow V\left[ \left[ t\right] \right] \\ 
x\bigotimes y\rightarrow \mu _t\left( x\bigotimes y\right)
=\sum_{m=0}^\infty \mu _m\left( x\bigotimes y\right) t^m
\end{array}
\ \text{with }\mu _m\in Hom\left( V\otimes V,V\right)
$$
$$
\Delta _t: 
\begin{array}{l}
V\rightarrow V\left[ \left[ t\right] \right] \otimes V\left[ \left[ t\right]
\right] =\left( V\otimes V\right) \left[ \left[ t\right] \right] \\ 
x\rightarrow \Delta _t\left( x\right) =\sum_{m=0}^\infty \Delta _m\left(
x\right) t^m
\end{array}
\ \text{with }\Delta _m\in Hom\left( V,V\otimes V\right)
$$

and they satisfy

\begin{itemize}
\item  $\mu _{t}$ is associative

\item  $\Delta _{t}$ is coassociative

\item  $(\mu _{t}\otimes \mu _{t})\circ (id\otimes \tau \otimes id)\circ\left( \Delta
_{t}\otimes \Delta _{t}\right) =\Delta _{t}\circ \mu _{t}\quad $where $\tau $ is
the twist map.
\end{itemize}

Two deformations $H_{t}$ and $H'_{t}$ of $H$ are said equivalent if there exits a formal
isomorphism $\phi _{t}=\phi _{0}+t\phi _{1}+...+t^{k}\phi_ {k}+...$ where 
$\phi _{0}=Id_{V}$ and $\phi_ {k} \in End(V)$ such that $H'_{t}=\phi _{t} \cdot H_{t}$ (as
defined in section $2$.). A deformation $H_{t}$ of $H$ is said trivial if $H_{t}$ is equivalent to $H$.
 
We give here some helpful remarks on deformations of Hopf algebras.

\begin{remark}

\begin{itemize}
\item If $H_{t}$ is a deformation of a Hopf algebra $H$ then the dual of $H_{t}$ is a
deformation of the dual of $H$. This follows from the linearity of the operations.
\item The deformation of pointed algebra in not always a pointed Hopf algebra.
Consider in $K^{p^3}$ where $p$ is prime and odd number the Hopf algebra defined by $A_t =
\frac{K \left\langle x, y, g \right\rangle}{ \left\langle g^{p} - 1, x^p, y^p, 
gx - qxg, gy - q^{-1} yg, xy - q^{- 1} yx - t ( 1 - g^2 ) \right\rangle}$,
 $\Delta g = g\otimes g$, $\Delta x = x \otimes 1 + g \otimes x$, 
 $\Delta y = y \otimes 1 + g\otimes y$, where t is a parameter in $K$ and $q$ is 
 a primitive $p$-th root of unity. We set $B_t = A_t^{\star}$ (the dual Hopf algebra) and $B_0 = A_0^{\star}$.
Since $A_t$ is a deformation of $A_0$ then $B_t$ is a deformation of $B_0$.
The Hopf algebra $B_0$ is pointed while its deformation $B_t$ is not pointed.
\item The number of grouplike elements of Hopf algebra may decreases by deformation. (see
the previous example)
\item The number of primitive idempotents of Hopf algebra do not decrease by deformation
(see\cite{Goze-makhlouf}). More generally, the order of a, element do not decrease by deformation.
\end{itemize}
\end{remark}

The following proposition gives a connection between degeneration and
deformation.

\begin{proposition}
\textit{If} $H_{0}$ \textit{is a degeneration of} $H_{1}$ \textit{then} $H_{1}$
 \textit{is a deformation of} $H_{0}.$
\end{proposition}

In fact, let $H_{0}=$ $\lim_{t\rightarrow 0}f_{t}\cdot H$ be a degeneration
of $H$ then $H_{t}=f_{t}\cdot H$ is a deformation of $H_{0}$.

\begin{remark}
\textbf{\ }The converse is, in general, false. We take from \cite{Beattie} 
\cite{Andrus-Schneider} the following example : for a fixed primitive $p$-th root of unity $\lambda $ and $t\in K^{*}$, the
family of Hopf algebra $H_{t}$ generated by $c,x_{1},x_{2}$ and subject to
\begin{center}
\begin{eqnarray*}
c^{p^{2}}=1,x_{1}^{p}=c^{p}-1,x_{2}^{p}=c^{p}-1, \\
x_{1}c =\lambda ^{-1}cx_{1},x_{2}c=\lambda cx_{2},x_{2}x_{1}=\lambda
x_{1}x_{2}+t(c^{2}-1) \\
\Delta (c) =c\otimes c,\Delta (x_{1})=c\otimes x_{1}+x_{1}\otimes 1,\Delta
(x_{2})=c\otimes x_{2}+x_{2}\otimes 1,
\end{eqnarray*}
\end{center}
The Hopf algebras $H_{t}$ and $H_{s}$ are isomorphic if and only if $ts^{-1}$
is a $p$-th root of unity. The Hopf algebra $H_{t}$ exists when $t\neq 0$
but not its limit when $t$ tends to $0$.   
\end{remark}

\subsection{Rigid Hopf algebras and irreducible components}

\begin{definition}
A Hopf algebra $H$ is said rigid if and only if every deformation of $H$ is
trivial.
\end{definition}

\begin{remark}
\end{remark}

\begin{itemize}
\item  The definition may be rephrased by ''The orbit of $H$ is Zariski
open''.

\item  If the second cohomological group of the Hopf algebra $H$ is trivial then $H$ is rigid.

\item  Every semisimple Hopf algebra is rigid because its second
cohomological groups is trivial \cite{Gerstehaber92}. 
Thus the group Hopf algebras are rigid. 

\item  The Zariski open orbits have a special interest in the geometric
study of $Hopf_{n}$, Their closure determines an irreducible component.

\item  Two non isomorphic rigid Hopf algebras belong to different
irreducible components.

\item  There is a finite number of open orbits because the algebraic variety 
$Hopf_{n}$ decomposes in a finite number of irreducible components.
\end{itemize}

\begin{theorem}

The dual $H^{\star}$, of a finite dimensional rigid Hopf algebra $H$,
is rigid.
\end{theorem}
\begin{proof}

Suppose that $H$ is rigid and $H^{\star}$ is not rigid. This means
that it exists a deformation $H^{\star}_{t}$ such that $H^{\star}_{t}$ is not isomorphic to $H^{*}$.\\
Since $(H^{\star}_{t})^{\star}$ is a deformation of $H^{\star \star}$, which is 
isomorphic to $H$ and $H$ rigid then
$(H^{\star}_{t})^{\star}
\simeq H$. By duality $(H^{\star}_{t})^{\star \star} \simeq H^{\star}$, which is equivalent 
to $(H_{t})^{\star} \simeq H^{\star} $, contradicting the hypothesis. 
\end{proof}

\begin{corollary}
The rigidity of a finite-dimensional Hopf algebra is equivalent to the rigidity of its dual.
\end{corollary}

\begin{corollary}
If a Hopf algebra $H_{0}$ is a degeneration of $H$ then $H^{\star}_{0}$ is a 
degeneration of $H^{\star}$.
\end{corollary}

\begin{proof}
Since $H_{0}$ is a degeneration of $H$ then it exists a deformation $H_{t}$ of $H$ which
tends to $H_{0}$ when $t$ tends to $0$. By duality $H^{\star}_{t}$ is a deformation of
$H^{\star}$ and by the linearity of the operations the limit of $H^{\star}_{t}$ is $H^{\star}_{0}$. 
\end{proof}

The following proposition give a necessary condition for the rigidity of of graded Hopf algebra.

\begin{proposition}
Let $H$ be a filtred Hopf algebra and $gr(H)$ its associated graded Hopf algebra.
If $H$ is not isomorphic to $gr(H)$ then $gr(H)$ is not rigid.  
\end{proposition}
\begin{proof}
Since $gr(H)$ is a degeneration of $H$ (theorem (3.3.1)) then $H$ is a deformation of $gr(H)$
(proposition 4.1.2). 
Thus, the graded Hopf algebra is not rigid if $H$ is not isomorphic to $gr(H)$. 
\end{proof}
\section{Rigidity and irreducible components in $Hopf_{n}$}
The general classification, up to isomorphism, of Hopf algebras is not
known. However, the complete classification is known for dimension $n$, $n<14
$ and for dimension $p$ and $p^{2}$ where $p$ is prime, see \cite{williams} 
\cite{stefan} \cite{Masuoka} \cite{fukuda}\cite{natale}\cite{Zhu1} and \cite{Ng}.
If the dimension is any prime number $p$ there is only one Hopf algebra, the Hopf group algebra
$KZ_{p}$. Then, $Hopf_{p}$ is formed by a unique irreducible component given by a
Zariski open orbit. In the following , we consider the varieties  $Hopf_{p^{2}}$, where $p$ is prime, and 
 $Hopf_{n}$, where $n<14$.
\subsection{Rigidity and irreducible components in $Hopf_{p^{2}}$, $p$ prime}
It is known from \cite{Ng} that in $Hopf_{p^{2}}$, where $p$ is any prime number, every
Hopf algebra is isomorphic to one of the following Hopf algebras :
\begin{enumerate}
\item $K[Z_{p^{2}}]$
\item $K[Z_{p}] \times K[Z_{p}]$
\item $T_{p^2}$ a Taft Hopf algebra, which is defined by

$$
\frac{K\left\langle x,y\right\rangle }{\left\langle x^p,\ y^p-1,\
xy-qyx\right\rangle }
$$

where $q$ is the primitive root of unity of order p.

The coalgebra algebra structure and the antipode are determined by

$\qquad \Delta (y)=y\otimes y,\quad \Delta (x)=x\otimes y+1\otimes x,$

$\qquad \varepsilon (x)=0,\quad \varepsilon (y)=1.$

$\qquad S(x)=-xy^{-1},\quad S(y)=y^{-1}.$
\end{enumerate}

\begin{theorem}
Every Hopf algebra $H$ in $Hopf_{p^{2}}$, where $p$ is any prime number, is rigid.

Therefore, the algebraic variety $Hopf_{p^{2}}$ is an union of $p+1$ Zariski open orbits.
\end{theorem}
\begin{proof}
The Hopf algebra $K[Z_{p^{2}}]$ and $K[Z_{p}] \times K[Z_{p}]$ are semisimple then
rigid. A Taft Hopf algebra cannot be deformed in a commutative Hopf algebra, then also rigid.

The number of irreducible components corresponds to the two semisimple Hopf algebras and
the $p-1$ non isomorphic Taft Hopf algebras.
\end{proof}
\subsection{Rigidity and irreducible components in $Hopf_{n}$ $n<14$}

For $n<12$, the Hopf algebra was classified by Williams and
reconsidered by Stefan. For $n=12$, the semisimple case  was studied by Fukuda and the
classification was completed by Natale. In the following, we recall the classification
 (see for the details \cite{stefan}, \cite{fukuda}, \cite{natale}).

 Let $Z_n$ denotes the cyclic group, $D_n$ the dihedral group, $S_n$ the
symmetric group, $H_4$ the quaternion group and $Al$ the alternate group.
Let $KG$ be the Hopf algebra of the group $G$ and $\left( KG\right) ^{*}$
its dual. Let $T_{n}$ be a Taft Hopf algebra described above.

\begin{theorem}
If $H$ is Hopf algebra of dimension $n<14$, then $H$ is isomorphic with
one and only one of the following Hopf algebra
\end{theorem}

\begin{itemize}
\item  
$\mathbf{n\in \{2,3,5,7,11,13}\}$
The group Hopf algebra $KZ_{n}$.

\item  $\mathbf{n=4.}$ 
The semisimple Hopf algebras $KZ_{4}$ and $%
K\left( Z_{2}\times Z_{2}\right) ,$ and the Taft-Sweedler Hopf algebra $T_{4}.$

\item  $\mathbf{n=6.}$
$KZ_6$, $KS_3$ and $\left( KS_3\right) ^{*}$

\item  $\mathbf{n=8}$. 
\textit{The semisimple Hopf algebras }\textbf{\ :} $K(Z_{2}\times
Z_{2}\times Z_{2})$, $K(Z_{2}\times Z_{4})$, $KZ_{8}$, $KD_{4}$, $\left(
KD_{4}\right) ^{*}$, $KH_{4}$, $\left( KH_{4}\right) ^{*}$ and $A_{8}$.
\textit{The nonsemisimple Hopf algebras  :}
$A_{C_2}$, 
$, 
A_{C_4}^{\prime }
$, 
$
A_{C_4}^{\prime \prime }
$, 
$
A_{C_4,q}^{\prime \prime \prime }
$
(where $q$ is the primitive root of unity of order 4),
$\left( A_{C_4}^{\prime \prime }\right) ^{*}$, and 
$
A_{C_2\times C_2}
$

\item  $\mathbf{n=9.}$
$KZ_9$, $K(Z_3\times Z_3)$ and the Taft Hopf algebras $T_9$.

\item  $\mathbf{n=10.}$
$KZ_{10}$, $KD_5$ and $\left( KD_5\right) ^{*}.$

\item  $\mathbf{n=12.}$
\textit{The semisimple Hopf algebras  : }$KZ_{12},$ $K(Z_{6}\times Z_{2})$%
, $K(Z_{4}\times Z_{3})$, $KD_{6}$, $\left( KD_{6}\right) ^{*}$, $Al_{4}$, $%
\left( Al_{4}\right) ^{*}$, $A_{+}$ and $A_{-}$ 

\textit{The nonsemisimple Hopf algebras  :} 
$
A_0
$, $
A_1
$
$
B_0
$
$
B_1
$  and 
$A^{\star}_{1}$.
\end{itemize}

\begin{theorem}
\textit{Every Hopf algebras of dimension }$n<14$ \textit{ is rigid except }
$A_{c_{4}}^{^{\prime }}$ \textit{in }
$Hopf_{8}$, 
\textit{ and } $A_{0}$\textit{ , } $A^{\star}_{0}\simeq B_{1}$ \textit{in} $Hopf_{12}$, .\newline
\textit{Furthermore,}\textit{\ in }$Hopf_{8}$, $A_{c_{4}}^{^{\prime }}$
\textit{\ is a degeneration of }$A_{c_{4}}^{^{\prime \prime }}$\textit{\ and in} 
$Hopf_{12}$, $A_{0}$\textit{\ is a
degeneration of }$A_{1} $\textit{and } $B_{1}$\textit{\ is a
degeneration of }$A^{\star}_{1}$\textit{. }
\end{theorem}

\begin{proof}
Following the previous remark, every semisimple Hopf algebra is rigid.

For $n\in \{$2,3,5,7,11,13$\},$ The Hopf algebras are all group Hopf algebras,
then rigid.

For $n=4$ and $n=9$, see theorem (5.1).

For $n=6,$ all Hopf algebras are semisimple, then rigid.

For $n=8,$ the semisimple Hopf algebras are rigid. 

The family $A_{c_{4},t}^{^{\prime }}=\frac{K\left\langle x,g\right\rangle }
{\left \langle g^{4}-1,x^{2}+t(1-g^{2}),gx+xg\right\rangle }$ 
is a deformation of $
A_{c_{4}}^{^{\prime }}$ and is isomorphic when $t\neq 0$ to 
$A_{c_{4}}^{^{\prime \prime }}$. Then $A_{c_{4}}^{^{\prime }}$ is a
degeneration of $A_{c_{4}}^{^{\prime \prime }}.$ The order of $g$ and the 
$(a,b)$-primitivity of $x$ don't allow deformation of other pointed
nonsemisimple Hopf algebras. 

For $n=10,$ all the algebras are semisimple, then rigid.

For $n=12,$ the semisimple Hopf algebras are rigid. 

The family $A_{t}=%
\frac{K\left\langle x,g\right\rangle }{\left\langle
g^{6}-1,x^{2}+t(g^{2}-1),gx+xg\right\rangle }$ is isomorphic to $A_{1}$ when 
$t\neq 0,$  and tends to $A_{0}$ when $t$ tends to $0$. then $A_{0}$ is a 
degeneration of $A_{1}$.

By duality $A^{\star}_{t}$ is a deformation of $A^{\star}_{0}$. Then $A^{\star}_{1} \simeq
A^{\star}_{t}$ ($t \neq 0$) is a deformation of $B_{1}$ ($B_{1} \simeq A^{\star}_{0}$) and
$B_{1}$ is a degeneration of $A^{\star}_{1}$. therefore $A_{0}$ and $B_{1}$ are not rigid
and all the others are rigid.
\end{proof}

\begin{corollary}
The following table gives the number of irreducible components of $Hopf_{n}$
for $n<14$

$$
\begin{tabular}{|c|c|}
\hline
$\text{dimension}$ & $\text{number of irreducible components of }Hopf_{n}$
\\ \hline
$n\in \{2,3,5,7,11,13\}$ & $1$ \\ \hline
$n=4$ & $3$ \\ \hline
$n=6$ & $3$ \\ \hline
$n=8$ & $14$ \\ \hline
$n=9$ & $4$ \\ \hline
$n=10$ & $3$ \\ \hline
$n=12$ & $14$ \\ \hline
\end{tabular}
$$
\end{corollary}
\begin{remark}
\begin{itemize}
\item The orbit dimension's of the $n$-dimensional semisimple
Hopf algebra is $n^{2}$.
\item It is interesting to see whether, for the nonsemisimple rigid Hopf algebras,
the second cohomological group is trivial.
\end{itemize}
\end{remark}
\section{Degenerations with $f_t=v+t\cdot w$}

The aim of this section is to find a necessary and sufficient
conditions such that a degeneration of a given Hopf algebra $H=\left( V,\mu
,\Delta ,\eta ,\varepsilon ,S\right) $ exists$.$ Let $f_{t}=v+tw$ be a
family of linear maps where $v$ is a singular linear map, $w$ is a regular
linear map and $t$ is a parameter in $K$. We suppose $v$ singular because when $v$ is regular, the family $%
f_{t}$ corresponds to isomorphisms. We can also set $w=id$ because $%
f_{t}=v+tw=\left( v\circ w^{-1}+t\right) \circ w$ which is isomorphic to $%
v\circ w^{-1}+t$. Then with no loss of generality we  consider the family $%
f_{t}=\varphi +t\cdot id$ from $V$ into $V$ where $\varphi $ is a singular
map and $t$ is in open set containing $0$. The vector space $V$ can be
decomposed by $\varphi $ under the form $V_{R}\oplus V_{N}$ where $V_{R}$
and $V_{N}$ are $\varphi $-invariant defined in a canonical way such that $%
\varphi $ is surjective on $V_{R}$ and nilpotent on $V_{N}$. Let $q$ be the
smallest integer such that $\varphi ^{q}\left( V_{N}\right) =0$. The inverse
of $f_{t}$ exists on $V_{R}$ and is equal to $\varphi ^{-1}\left( t\varphi
^{-1}+id\right) ^{-1}$. But on $V_{N},$ since $\varphi ^{q}=0$, it is given
by 
$$
\frac{1}{\varphi +t\cdot id}=\frac{1}{t}\cdot \frac{1}{\varphi /t+id}=\frac{1%
}{t}\cdot \sum_{i=0}^{\infty }\left( -\frac{\varphi }{t}\right) ^{i}=\frac{1%
}{t}\cdot \sum_{i=0}^{q-1}\left( -\frac{\varphi }{t}\right) ^{i}
$$

It follows :

\smallskip \textbf{Fitting lemma.} \textit{Let }$f_{t}=\varphi +t\cdot id$%
\textit{\ be a family of linear maps from }$V$\textit{\ into }$V$\textit{\
where }$\varphi $\textit{\ is a singular map. Then, }$V=V_{R}\oplus V_{N}$%
\textit{\ where }$V_{R}$\textit{\ and }$V_{N}$\textit{\ are }$\varphi $%
\textit{-invariant and }$\varphi $\textit{\ is surjective on }$V_{R}$\textit{%
\ and nilpotent on }$V_{N}$\textit{. The inverse of }$f_{t}$\textit{\ is
defined by}

$$
f_{t}^{-1}=\left\{ 
\begin{array}{c}
\varphi ^{-1}\left( t\varphi ^{-1}+id\right) ^{-1}\quad \text{on }V_{R} \\ 
\frac{1}{t}\cdot \sum_{i=0}^{q-1}\left( -\frac{\varphi }{t}\right) ^{i}\quad 
\text{on }V_{N}
\end{array}
\right.
$$
\textit{where }$q$\textit{\ is the smallest integer such that }$\varphi
^{q}\left( V_{N}\right) =0.$

\subsection{Degeneration of an algebra}

Let $f_t=\varphi +t\cdot id$ be a family of linear maps on $V$, where $%
\varphi $ is a singular map. The action of $f_t$ on $\mu $ is defined by $%
f_t\cdot \mu =f_t^{-1}\circ \mu \circ f_t\otimes f_t$ then 
$$
\begin{array}{c}
f_t\cdot \mu \left( x\otimes y\right) =f_t^{-1}\circ \mu \left( f_t\left(
x\right) \otimes f_t\left( y\right) \right) \\ 
=f_t^{-1}\left( \mu \left( \varphi \left( x\right) \otimes \varphi \left(
y\right) \right) +t\left( \mu \left( \varphi \left( x\right) \otimes
y\right) +\mu \left( x\otimes \varphi \left( y\right) \right) \right)
+t^2\mu \left( x\otimes y\right) \right)
\end{array}
$$
Since every element $v$ of $V$ decomposes in $v=v_R+v_N$, we set

$$
\begin{array}{c}
A=\mu \left( x\otimes y\right) =A_R+A_N, \\ 
B=\mu \left( \varphi \left( x\right) \otimes y\right) +\mu \left( x\otimes
\varphi \left( y\right) \right) =B_R+B_N \\ 
C=\mu \left( \varphi \left( x\right) \otimes \varphi \left( y\right) \right)
=C_R+C_N
\end{array}
$$

Then 
$$
f_t\cdot \mu \left( x\otimes y\right) =\varphi ^{-1}\left( t\varphi
^{-1}+id\right) ^{-1}\left( t^2A_R+tB_R+C_R\right) +\frac 1t\cdot
\sum_{i=0}^{q-1}\left( -\frac \varphi t\right) ^i\left(
t^2A_N+tB_N+C_N\right)
$$

If $t$ goes to $0,$ then $\varphi ^{-1}\left( t\varphi ^{-1}+id\right)
^{-1}\left( t^2A_R+tB_R+C_R\right) $ goes to $\varphi ^{-1}\left( C_R\right)
.$ The limit of the second term is :

\begin{center}
$\lim _{t\rightarrow 0}\frac 1t\cdot \sum_{i=0}^{q-1}\left( -\frac \varphi
t\right) ^i\left( t^2A_N+tB_N+C_N\right) $

$=\lim _{t\rightarrow 0}\left( \frac{id}t-\frac \varphi {t^2}+\frac{\varphi
^2}{t^3}-\cdots \left( -1\right) ^{q-1}\frac{\varphi ^{q-1}}{t^q}\right)
\left( t^2A_N+tB_N+C_N\right) $

$=\lim _{t\rightarrow 0}tA_N-\varphi \left( A_N\right) +\left( \frac{id}%
t-\frac \varphi {t^2}+\frac{\varphi ^2}{t^3}-\cdots \left( -1\right) ^{q-1}%
\frac{\varphi ^{q-3}}{t^{q-2}}\right) \varphi ^2\left( A_N\right) +$

$B_N-\left( \frac{id}t-\frac \varphi {t^2}+\frac{\varphi ^2}{t^3}-\cdots
\left( -1\right) ^{q-1}\frac{\varphi ^{q-2}}{t^{q-1}}\right) \varphi \left(
B_N\right) +$

$\left( \frac{id}t-\frac \varphi {t^2}+\frac{\varphi ^2}{t^3}-\cdots \left(
-1\right) ^{q-1}\frac{\varphi ^{q-1}}{t^q}\right) C_N$

$=\lim _{t\rightarrow 0}tA_N+B_N-\varphi \left( A_N\right) +\qquad \qquad
\qquad $

$\left( \frac{id}t-\frac \varphi {t^2}+\frac{\varphi ^2}{t^3}-\cdots \left(
-1\right) ^{q-1}\frac{\varphi ^{q-1}}{t^q}\right) \left( \varphi ^2\left(
A_N\right) -\varphi \left( B_N\right) +C_N\right) $
\end{center}

This limit exists if and only if 
$$
\varphi ^2\left( A_N\right) -\varphi \left( B_N\right) +C_N=0
$$
which is equivalent to 
$$
\left( 11\right) \quad \varphi ^2\left( \mu \left( x\otimes y\right)
_N\right) -\varphi \left( \mu \left( \varphi \left( x\right) \otimes
y\right) _N\right) -\varphi \left( \mu \left( x\otimes \varphi \left(
y\right) \right) _N\right) +\mu \left( \varphi \left( x\right) \otimes
\varphi \left( y\right) \right) _N=0
$$
And the limit is $B_N-\varphi \left( A_N\right) .$

\begin{proposition}
\textbf{\ }\textit{The degeneration of the algebra }$\mu $\textit{\ exists
if and only if the condition 
$$
\left( 11\right) \quad \varphi ^{2}\circ \mu _{N}-\varphi \circ \mu
_{N}\circ \varphi \otimes id-\varphi \circ \mu _{N}\circ id\otimes \varphi
+\mu _{N}\circ \varphi \otimes \varphi =0
$$
\ where }$\mu _{N}\left( x,y\right) =\left( \mu \left( x,y\right) \right)
_{N}$, \textit{holds. And it is defined by} 
$$
\mu _{0}=\varphi ^{-1}\circ \mu _{R}\circ \varphi \otimes \varphi +\mu
_{N}\circ \varphi \otimes id+\mu _{N}\circ id\otimes \varphi -\varphi \circ
\mu _{N}
$$
\end{proposition}

\subsection{Degeneration of a coalgebra}

Let $f_t=\varphi +t\cdot id$ be a family of linear maps on $V$, where $%
\varphi $ is a singular map. The action of $f_t$ on $\Delta $ is defined by $%
f_t\cdot \Delta =f_t^{-1}\otimes f_t^{-1}\circ \Delta \circ f_t$ then 
$$
f_t\cdot \Delta \left( x\right) =f_t^{-1}\otimes f_t^{-1}\circ \Delta \left(
f_t\left( x\right) \right) =t\cdot f_t^{-1}\otimes f_t^{-1}\circ \Delta
\left( x\right) +f_t^{-1}\otimes f_t^{-1}\left( \Delta \left( \varphi \left(
x\right) \right) \right)
$$
Setting $\Delta \left( x\right) =x^{\left( 1\right) }\otimes x^{\left(
2\right) }$, $\Delta \left( \varphi \left( x\right) \right) =\varphi \left(
x\right) ^{\left( 1\right) }\otimes \varphi \left( x\right) ^{\left(
2\right) }$, \newline
and for $i=1,2\quad x^{\left( i\right) }=x_R^{\left( i\right) }+x_N^{\left(
i\right) }$,$\ \varphi \left( x\right) ^{\left( i\right) }=\varphi \left(
x\right) _R^{\left( i\right) }\otimes \varphi \left( x\right) _N^{\left(
i\right) }.$\newline
Then \newline

$f_t\cdot \Delta \left( x\right) =t\cdot \left( \left( f_t^{-1}\left(
x_R^{\left( 1\right) }\right) +f_t^{-1}\left( x_N^{\left( 1\right) }\right)
\right) \otimes \left( f_t^{-1}\left( x_R^{\left( 2\right) }\right)
+f_t^{-1}\left( x_N^{\left( 2\right) }\right) \right) \right) +$

+$\left( f_t^{-1}\left( \varphi \left( x\right) _R^{\left( 1\right) }\right)
+f_t^{-1}\left( \varphi \left( x\right) _N^{\left( 1\right) }\right) \right)
\otimes \left( f_t^{-1}\left( \varphi \left( x\right) _R^{\left( 2\right)
}\right) +f_t^{-1}\left( \varphi \left( x\right) _N^{\left( 2\right)
}\right) \right) $\newline
Setting $\psi =\varphi ^{-1}\left( t\varphi ^{-1}+id\right) ^{-1}$, then
when $t\rightarrow 0,$ $\psi =\varphi ^{-1}$ and\newline
$f_t\cdot \Delta \left( x\right) =t\cdot \psi \left( x_R^{\left( 1\right)
}\right) \otimes \psi \left( x_R^{\left( 2\right) }\right) +\psi \left(
\varphi \left( x\right) _R^{\left( 1\right) }\right) \otimes \psi \left(
\varphi \left( x\right) _R^{\left( 2\right) }\right) +$

$+\sum_{i=0}^{q-1}\psi \left( x_R^{\left( 1\right) }\right) \otimes \left(
-\frac \varphi t\right) ^i\left( x_N^{\left( 2\right) }\right) +\left(
-\frac \varphi t\right) ^i\left( x_N^{\left( 1\right) }\right) \otimes \psi
\left( x_R^{\left( 2\right) }\right) +$

$+\frac 1t\sum_{i=0}^{q-1}\psi \left( \varphi \left( x\right) _R^{\left(
1\right) }\right) \otimes \left( -\frac \varphi t\right) ^i\left( \varphi
\left( x\right) _N^{\left( 2\right) }\right) +\left( -\frac \varphi t\right)
^i\left( \varphi \left( x\right) _N^{\left( 1\right) }\right) \otimes \psi
\left( \varphi \left( x\right) _R^{\left( 2\right) }\right) +$

+$\frac 1t\sum_{i=0}^{q-1}\sum_{j=0}^{q-1}\left( -\frac \varphi t\right)
^j\left( x_N^{\left( 1\right) }\right) \otimes \left( -\frac \varphi
t\right) ^i\left( x_N^{\left( 2\right) }\right) $+$\frac
1{t^2}\sum_{i=0}^{q-1}\sum_{j=0}^{q-1}\left( -\frac \varphi t\right)
^j\left( \varphi \left( x\right) _N^{\left( 1\right) }\right) \otimes \left(
-\frac \varphi t\right) ^i\left( \varphi \left( x\right) _N^{\left( 2\right)
}\right) $

Then\newline
$f_t\cdot \Delta \left( x\right) =\psi \otimes \psi \left( t\cdot
x_R^{\left( 1\right) }\otimes x_R^{\left( 2\right) }+\varphi \left( x\right)
_R^{\left( 1\right) }\otimes \varphi \left( x\right) _R^{\left( 2\right)
}\right) $

$+\sum_{i=0}^{q-1}\frac 1{\left( -t\right) ^i}\left( \psi \otimes \varphi
^i\left( x_R^{\left( 1\right) }\otimes x_N^{\left( 2\right) }\right)
+\varphi ^i\otimes \psi \left( x_N^{\left( 1\right) }\otimes x_R^{\left(
2\right) }\right) \right) +$

$-\sum_{i=0}^{q-1}\frac 1{\left( -t\right) ^{i+1}}\left( \psi \otimes
\varphi ^i\left( \varphi \left( x\right) _R^{\left( 1\right) }\otimes
\varphi \left( x\right) _N^{\left( 2\right) }\right) +\varphi ^i\otimes \psi
\left( \varphi \left( x\right) _N^{\left( 1\right) }\otimes \varphi \left(
x\right) _R^{\left( 2\right) }\right) \right) +$

+$\sum_{i=0}^{q-1}\sum_{j=0}^{q-1}\frac{-1}{\left( -t\right) ^{i+j+1}}%
\varphi ^j\otimes \varphi ^i\left( x_N^{\left( 1\right) }\otimes x_N^{\left(
2\right) }\right) $+$\sum_{i=0}^{q-1}\sum_{j=0}^{q-1}\frac{-1}{\left(
-t\right) ^{i+j+2}}\varphi ^j\otimes \varphi ^i\left( \varphi \left(
x\right) _N^{\left( 1\right) }\otimes \varphi \left( x\right) _N^{\left(
2\right) }\right) $\newline
By rewriting the sums one obtain\newline
$f_t\cdot \Delta \left( x\right) =\psi \otimes \psi \left( t\cdot
x_R^{\left( 1\right) }\otimes x_R^{\left( 2\right) }+\varphi \left( x\right)
_R^{\left( 1\right) }\otimes \varphi \left( x\right) _R^{\left( 2\right)
}\right) +$

$\psi \otimes id\left( x_R^{\left( 1\right) }\otimes x_N^{\left( 2\right)
}\right) +id\otimes \psi \left( x_N^{\left( 1\right) }\otimes x_R^{\left(
2\right) }\right) +$

$\sum_{i=0}^{q-1}\frac 1{\left( -t\right) ^{i+1}}\psi \otimes \varphi
^i\left( x_R^{\left( 1\right) }\otimes \varphi \left( x_N^{\left( 2\right)
}\right) -\varphi \left( x\right) _R^1\otimes \varphi \left( x\right)
_N^{\left( 2\right) }\right) +$

$\sum_{i=0}^{q-1}\frac 1{\left( -t\right) ^{i+1}}\varphi ^i\otimes \psi
\left( \varphi \left( x_N^{\left( 1\right) }\right) \otimes x_R^{\left(
2\right) }-\varphi \left( x\right) _N^{\left( 1\right) }\otimes \varphi
\left( x\right) _R^{\left( 2\right) }\right) +$

$\sum_{i=0}^{q-1}\frac 1{\left( -t\right) ^{i+1}}\varphi ^i\otimes id\left(
x_N^{\left( 1\right) }\otimes x_N^{\left( 2\right) }\right) +$

$\sum_{i=0}^{q-1}\sum_{j=0}^{q-1}\frac{-1}{\left( -t\right) ^{i+j+2}}\varphi
^j\otimes \varphi ^i\left( \varphi \left( x\right) _N^{\left( 1\right)
}\otimes \varphi \left( x\right) _N^{\left( 2\right) }-x_N^{\left( 1\right)
}\otimes \varphi \left( x_N^{\left( 2\right) }\right) \right) $\newline
This limit of $f_t\cdot \Delta \left( x\right) $ exists if and only if 
$$
\left( 12\right) \quad \left\{ 
\begin{array}{l}
x_R^{\left( 1\right) }\otimes \varphi \left( x_N^{\left( 2\right) }\right)
-\varphi \left( x\right) _R^1\otimes \varphi \left( x\right) _N^{\left(
2\right) }=0 \\ 
\varphi \left( x_N^{\left( 1\right) }\right) \otimes x_R^{\left( 2\right)
}-\varphi \left( x\right) _N^{\left( 1\right) }\otimes \varphi \left(
x\right) _R^{\left( 2\right) }=0 \\ 
x_N^{\left( 1\right) }\otimes x_N^{\left( 2\right) }=0 \\ 
\varphi \left( x\right) _N^{\left( 1\right) }\otimes \varphi \left( x\right)
_N^{\left( 2\right) }-x_N^{\left( 1\right) }\otimes \varphi \left(
x_N^{\left( 2\right) }\right) =0
\end{array}
\right. \qquad \forall x\in V
$$

\begin{proposition}
\textit{The degeneration of the coalgebra }$\Delta $\textit{\ exists if and
only if the condition }$\left( 12\right) $\textit{\ holds and it is defined
for all }$x\in V$ \textit{by} 
$$
\Delta _{0}\left( x\right) =\varphi ^{-1}\left( \varphi \left( x\right)
_{R}^{\left( 1\right) }\right) \otimes \varphi ^{-1}\left( \varphi \left(
x\right) _{R}^{\left( 2\right) }\right) +\varphi ^{-1}\left( x_{R}^{\left(
1\right) }\right) \otimes x_{N}^{\left( 2\right) }+x_{N}^{\left( 1\right)
}\otimes \varphi ^{-1}\left( x_{R}^{\left( 2\right) }\right) 
$$
\end{proposition}

\subsection{Degeneration of Hopf algebra}

Let $H=(V,\mu ,\eta ,\Delta ,\varepsilon ,S)$ be a Hopf algebra and $%
f_{t}=\varphi +t\cdot id$ be a family of linear maps of $V$, where $\varphi $
is a singular map. We suppose also that $\varphi $ decomposes the vector
space $V$ as $V=V_{R}+V_{N}$. Then the degeneration $H_{0}=\lim_{t%
\rightarrow 0}f_{t}\cdot H$ exists if and only if the conditions $\left(
11\right) $ $\left( 12\right) $ holds, the antipode remains for $H_{0}$. The
multiplication and the comultiplication defined by 
$$
\begin{array}{c}
\mu _{0}(x \otimes y)=\varphi ^{-1}\left( \mu \left( \varphi \left( x\right) \otimes
\varphi \left( y\right) \right) _{R}\right) +\mu \left( \varphi \left(
x\right) \otimes y\right) _{N}+\mu \left( x\otimes \varphi \left( y\right)
\right) _{N}-\varphi \left( \mu \left( x\otimes y\right) _{N}\right)  \\ 
\Delta _{0}\left( x\right) =\varphi ^{-1}\left( \varphi \left( x\right)
_{R}^{\left( 1\right) }\right) \otimes \varphi ^{-1}\left( \varphi \left(
x\right) _{R}^{\left( 2\right) }\right) +\varphi ^{-1}\left( x_{R}^{\left(
1\right) }\right) \otimes x_{N}^{\left( 2\right) }+x_{N}^{\left( 1\right)
}\otimes \varphi ^{-1}\left( x_{R}^{\left( 2\right) }\right) 
\end{array}
$$
satisfy the condition

($\mu _0\otimes \mu _0)\circ (id\otimes \tau \otimes id)\circ \left( \Delta _0\otimes
\Delta _0\right) =\Delta _0 \circ \mu _0\quad $where $\tau $ is the twist map

$
\mu_0 \circ (S\otimes Id)\circ \Delta_0 =\mu_0 \circ (Id\otimes S)\circ
\Delta_0 =\eta \circ \varepsilon
$

\end{document}